\documentclass[a4paper]{amsart}
\usepackage{amssymb}  
\usepackage{amsmath}  
\usepackage{amscd}
\usepackage{graphicx} 
\title[Regular homotopy classes of immersions
of $3$-manifolds into $5$-space]
{Regular homotopy classes of immersions\\
of $3$-manifolds into $5$-space}
\author{Osamu Saeki, Andr\'as Sz\H{u}cs and 
Masamichi TAKASE}
\thanks{The first author has been partially supported by
Grant-in-Aid for Scientific Research (No.~11440022),
Ministry of Education, Science and Culture, Japan.
The second author has been supported by the Grant
OTKA T029795.}
\address{Department of Mathematics, Graduate School of 
Science, Hiroshima University, Higashi-Hiroshima 739-8526,
Japan} 
\email{saeki@math.sci.hiroshima-u.ac.jp}
\address{Department of Analysis, ELTE, 
R\'ak\'oczi \'ut 5-7, Budapest 1088, Hungary} 
\email{szucs@cs.elte.hu}
\address{Graduate School of Mathematical Sciences, 
University of Tokyo, 3-8-1 Komaba, Meguro-ku, Tokyo, 
153-8914, Japan} 
\email{takase@ms.u-tokyo.ac.jp}
\newtheorem{theorem}{Theorem}[section]
\newtheorem{corollary}[theorem]{Corollary}
\newtheorem{lemma}[theorem]{Lemma}
\newtheorem{proposition}[theorem]{Proposition}
\theoremstyle{definition}
\newtheorem{definition}[theorem]{Definition}
\newtheorem{remark}[theorem]{Remark}

\newcommand{\punc}[1]{{#1}_{\circ}}
\newcommand{\sk}[2]{{\rm sk}_{#1}{#2}}

\newcommand{\R}{\mathbf{R}}
\newcommand{\Z}{\mathbf{Z}}
\newcommand{\Q}{\mathbf{Q}}
\newcommand{\Imm}{\mathop{\mathrm{Imm}}\nolimits}
\newcommand{\Emb}{\mathop{\mathrm{Emb}}\nolimits}
\newcommand{\Ker}{\mathop{\mathrm{Ker}}\nolimits}

\newcommand{\Int}{\mathop{\mathrm{int}}\nolimits}

\renewcommand{\tilde}{\widetilde}
\renewcommand{\hat}{\widehat}
\renewcommand{\setminus}{\smallsetminus}
\begin{document}
\begin{abstract}
We give geometric formulae which enable us
to detect (completely in some cases) the regular homotopy 
class of an immersion with trivial normal bundle of a 
closed oriented 
$3$-manifold into $5$-space.
These are analogues of the geometric formulae 
for the Smale invariants due to Ekholm and the second 
author. As a corollary,
we show that two embeddings into $5$-space of a closed 
oriented $3$-manifold with no $2$-torsion in the second 
cohomology are regularly homotopic if and only if they 
have Seifert surfaces with the same signature.
We also show that 
there exist two embeddings 
$F_0$ and $F_8 : T^3 \hookrightarrow \R^5$
of the $3$-torus $T^3$ with the following properties:
$(1)$ $F_0 \sharp h$ is regularly homotopic to $F_8$ for 
some immersion $h : S^3 \looparrowright \R^5$, and
$(2)$ the immersion $h$ as above cannot be chosen from a
regular homotopy class containing an embedding.
\end{abstract}
\maketitle
\section{Introduction}\label{sect:intro}
The Smale-Hirsch theory \cite{hirsch,smale} reduces the 
problem of 
enumerating regular homotopy classes of immersions
to homotopy theory. However, the problem of determining 
the regular homotopy class of a given immersion remains 
nontrivial. Here we 
address this problem in the case of immersions of 
$3$-manifolds into $\R^5$ having trivial normal bundles.
We concentrate on immersions with trivial normal bundles,
since we want to apply our results to the 
following question:
{\it Which regular homotopy classes of immersions of an 
oriented $3$-manifold into $\R^5$ contain embeddings?}
Note that embeddings of oriented $3$-manifolds 
into $\R^5$ have trivial normal bundles. 

\smallskip
\noindent{\it History.}\ \par
The problem of reading off the regular homotopy class of 
an immersion from its picture was posed by Smale 
\cite{smale}. For the particular
case of immersions of the $n$-sphere $S^n$ into $\R^{2n}$, 
it had already been shown 
by Whitney \cite{whitney} that the double points of an 
immersion $S^n\looparrowright\R^{2n}$ determine its regular 
homotopy class. This result has been extended by Ekholm 
\cite{ekholm2,ekholm3}
to immersions $S^n\looparrowright\R^{2n-1} (n \ge 4)$ and 
$S^n\looparrowright\R^{2n-2}.$

On the other hand, a result by Hughes and Melvin \cite{h-m} 
suggested that no further extension was possible. 
Namely they have shown that there are embeddings
$S^{4k-1}\hookrightarrow\R^{4k+1}$ not regularly homotopic 
to the standard one. Therefore the multiple points do not 
determine the regular homotopy class of an immersion.

Still the original problem posed by Smale has been answered 
by Ekholm and Sz\H ucs \cite{e-s} for immersions 
$S^{4k-1}\looparrowright\R^{4k+1}$ 
by considering ``singular Seifert surfaces'' bounded by the 
immersions. 

Immersions of 
closed oriented $3$-manifolds into $\R^5$ 
were investigated by Wu \cite{wu}. 
Using the  homotopy theoretical reduction provided by the 
Smale-Hirsch theory
he has shown that the regular homotopy classes of 
immersions of
a closed oriented $3$-manifold $M^3$ into $\R^5$
can be characterised by two invariants as follows.
\begin{itemize}
\item[(1)] The first one is a second cohomology class of 
$M^3$:
``half of the normal Euler class''.
\item[(2)] The second invariant is an integer 
(for immersions with trivial normal bundles).
\end{itemize}
The geometric meaning of these invariants was not very 
clear from Wu's work. 

\smallskip
\noindent{\it Now we can list our results:}
\renewcommand{\theenumi}{\roman{enumi}}
\renewcommand{\labelenumi}{\theenumi)}
\begin{enumerate}
\item
We shall give a geometric description of the first 
invariant
in the first part of Section~\ref{sect:C}.
We call it ``the Wu invariant''.
Note that for an embedding $M^3 \hookrightarrow \R^5$ the 
Wu invariant is
either zero or an element of order two in the second 
cohomology group
$H^2(M^3;\Z)$ of $M^3.$
\item 
We show that given any second cohomology class in 
$H^2(M^3;\Z)$, 
which is zero or has order two,
there exists an embedding 
$M^3\hookrightarrow\R^5$ of the oriented $3$-manifold 
$M^3$
with Wu invariant equal to the given cohomology class 
(Theorem~\ref{thm:C}).
\item  
We define an action of the group $\Imm[S^3,\R^5]$
of regular homotopy classes of immersions of $S^3$ into 
$\R^5$ 
on the set $\Imm[M^3,\R^5]_0$ of regular homotopy classes 
of immersions of the 
$3$-manifold $M^3$ into $\R^5$ with trivial normal bundles
by taking connected sums of immersions. 
In Section~\ref{sect:sharp}, we show that this action is 
effective
and its orbits coincide with the regular homotopy classes 
of immersions 
having the same Wu invariant.
\item\label{item:4}
In Section~\ref{sect:main}, we express the second invariant 
(which is an integer) through singularities of 
any singular Seifert surface. Actually we prove two 
formulae, analogous to those proved by Ekholm and Sz\H ucs 
\cite{e-s} for immersions of $S^3$ into $\R^5$.
\item 
In Section~\ref{sect:t3}, we prove a few rather surprising 
corollaries concerning the question mentioned above about 
the regular homotopy classes containing embeddings.
Although these corollaries hold more generally, we 
formulate them for the simplest manifold for which these 
surprising phenomena occur, namely for the 
$3$-dimensional torus $T^3$:
\vspace{2mm}

\noindent{\bf Corollary~\ref{cor:t3sz}.}
{\it There exist embeddings $F_0$ and 
$F_8 : T^3 \hookrightarrow \R^5$ of the $3$-torus $T^3$ 
with the following properties:
$(1)$ There is an
immersion $h : S^3 \looparrowright \R^5$ such that
$F_8$ is regularly homotopic to the connected sum $F_0 \sharp 
h$, but
$(2)$ this immersion $h$ cannot be chosen from the regular 
homotopy
class of an embedding.}\vspace{2mm}

\noindent{\bf Corollary~\ref{cor:t3}.}
{\it There exists an immersion $h:S^3\looparrowright\R^5$ 
with the following properties:
$(1)$ $h$ is not regularly homotopic to an embedding 
$S^3\hookrightarrow \R^5$, but
$(2)$ for any embedding $E:T^3\hookrightarrow\R^5$,
the connected sum $E \sharp h:T^3\looparrowright \R^5$
is again regularly homotopic to 
an embedding $T^3\hookrightarrow \R^5$.
}\vspace{2mm}

In other words, the double points of the immersion $h$ 
cannot be eliminated
by regular homotopy, 
but they can be eliminated after being
taken connected sum with any embedding of the $3$-torus.
\item
In Section~\ref{sect:cusp}, we give some corollaries 
concerning the 
number of cusps a singular Seifert surface can have.
\end{enumerate}

We will use the following notation throughout the paper.
By $\Imm[M,N]$ we denote the set of regular homotopy classes 
of immersions of a manifold $M$ into a manifold $N$
and by $\Emb[M,N]$ the subset of
$\Imm[M,N]$ consisting of all regular homotopy classes 
containing an
embedding. By $\Imm[M,N]_{\chi}$ we denote the subset of
$\Imm[M,N]$ consisting of all regular homotopy classes
of immersions with normal Euler class $\chi \in H^*(M;\Z)$.
In particular for an oriented $3$-manifold $M^3$ by 
$\Imm[M^3,\R^5]_0$ 
we denote the set of regular homotopy classes of immersions
$M^3 \looparrowright \R^5$
with trivial normal bundles.
By $M^3$ we will always denote a smooth, closed, oriented 
$3$-manifold.

We denote by $\Gamma_2(M^3)$ the finite subgroup of 
$H^2(M^3;\Z)$ 
formed by zero and all elements of order two.

The punctured manifold 
of $M^3$ will be denoted by $\punc M^3$; i.e.\ 
$\punc{M}^3:=M^3\setminus\Int D^3$
($D^3$ is a $3$-disk) 
and the $i$-skeleton of a manifold $M$ by $\sk{i}{M}$.

Throughout the paper, manifolds and maps are 
of class $C^\infty$.
The symbol ``$\approx$''
denotes an appropriate isomorphism between algebraic 
objects;
``$\simeq$'' and ``$\sim _r$'' mean
``homotopic'' and ``regularly homotopic'' respectively. 

\vspace{5mm}  
\section{Geometric formulae for the Smale invariants\\
given by Hughes-Melvin and by 
Ekholm-Sz\H{u}cs}\label{sect:rev}

Hughes \cite{hughes} has shown that
$\Imm[S^n,\R^N]$ has a group structure under connected sum 
and the Smale invariant 
$$\Omega:\Imm[S^n,\R^N] \to \pi_n(V_{N,n})$$  
actually gives a group isomorphism, where $V_{N,n}$
denotes the Stiefel manifold of $n$ frames in the 
$N$-space.

Hughes and Melvin \cite{h-m} investigated the problem: 
Which regular 
homotopy classes contain
embeddings in the set  
$\Imm[S^{4k-1},\R^{4k+1}]$? Their result for the case of 
embeddings of $S^3$ 
into $\R^5$ can be summarised by the following diagram:
$$\begin{array}{cccc}
\Omega: & \Imm[S^3, \R^5] & 
\stackrel{\approx}{\longrightarrow} 
&\pi _3 (V_{5,3})\approx \Z\\
 & \cup & & \cup \\
 & \Emb[S^3, \R^5] &\stackrel{\approx}{\longrightarrow} 
& 24\Z\\
 & f & \longmapsto & 3\sigma (V^4_f)/2,
\end{array}$$  
where $\Omega$ is the Smale invariant and
$\sigma (V^4_f)$ is the signature of an oriented Seifert 
surface 
$V^4_f$ for $f$. Here, ``$\approx$'' denotes a group 
isomorphism.

Their result has been extended by the third author to 
$\Z_2$-homology 3-spheres \cite{takase}.
Trying to extend it further was the starting point of 
the present paper.

As we mentioned in the introduction,

\begin{itemize}
\item[(1)]
this result shows the impossibility of determining the 
regular homotopy 
class of an immersion by its multiple points, and 
\item[(2)]
it has been shown in \cite{e-s} that Hughes-Melvin's 
formula stating 
that the Smale invariant of 
an embedding is equal to $3/2$ times the signature of 
any Seifert surface 
bounded by the embedding can be extended to immersions 
if we allow ``singular Seifert surfaces''.
\end{itemize}
Since we are going to generalise the result (2) to 
immersions of 
arbitrary oriented $3$-manifolds with trivial normal 
bundles, 
we recall it in detail.

Let $f : S^3 \looparrowright \R^5$ be an immersion and
$V^4$ an arbitrary compact oriented 4-manifold with
$\partial V^4 = S^3$.
The map $f$ extends to a
generic map $\tilde{f} : V^4 \to \R^5$ which has no 
singular points near the 
boundary. This map $\tilde{f}$ has isolated cusps, 
each of which has
a sign. Let us denote by
$\#\Sigma^{1,1}(\tilde{f})$ their algebraic number. 
Let us denote by $\Omega(f)$ the Smale invariant of $f$.
The first formula for the Smale invariant given by 
\cite{e-s} is the following:

\begin{theorem}[Theorem~1 (a) in \cite{e-s}]\label{thm2.1}
$$\Omega(f) = \frac{1}{2}(3\sigma(V^4)+
\#\Sigma^{1,1}(\tilde{f})).\quad\qed$$
\end{theorem}
\smallskip

Before giving the second formula for the Smale invariant
we recall some definitions from \cite{e-s}. 
Let $V^4$ be as above and let
$\hat{f} : V^4 \to \R^6_+$ be a generic map 
nonsingular near the boundary, 
such that
$\hat{f}^{-1}(\R^5) = \partial V^4 = S^3$ and
$\hat{f}|_{\partial V^4} \sim _r f$ in $\R^5$.
Here $\R^6_+$ is the upper half space in $\R^6$
and $\R^5 = \partial\R^6_+$.
%
The double point set $\Delta(\hat{f}) 
\subset V^4$ of the map $\hat f$ is 2-dimensional, and
on the boundary of $\hat{f}(\Delta(\hat{f}))$
lies the image $\Sigma$ of the singular points of 
$\hat{f}$, which
is 1-dimensional. Furthermore, triple points appear
as isolated points, each of which has a sign.
We denote the algebraic number of triple points of $\hat f$ 
by $t(\hat f)$.

\begin{definition}[%
Definition~4 in \cite{e-s}]\label{dfn:L}
Let $\Delta=\Delta(f)$ be the set of double points of $f$
and $f(\Delta)$ its image.
Then $f(\Delta)$ is a closed $1$-manifold in $\R^5$.
Since $S^3$ is $2$-connected, $f$ has a unique (up to 
homotopy) normal framing ($n_1,n_2$).
Let $f(\Delta)'$ be the manifold obtained from $f(\Delta)$ 
by shifting slightly along the vector field, 
$q\mapsto n_1(p_1)+n_1(p_2)$ for $q=f(p_1)=f(p_2)\in 
f(\Delta)$.
Define $L(f)$ to be the linking number in $\R^5$ of 
$f(S^3)$ and $f(\Delta)'$.
\end{definition}

\begin{definition}[%
Definition~5 in \cite{e-s}]\label{dfn:l}
Let $\Sigma'$ be a copy of $\Sigma$ 
shifted slightly along the outward normal vector field 
of $\Sigma$ in (the closure of) 
${\hat{f}(\Delta(\hat{f}))}$.
Define the linking number $l(\hat{f})\in \Z$ 
of the map $\hat{f}$
to be the linking number of $\Sigma'$ 
and $\hat{f}(V^4)$ in $\R^6_+$.
\end{definition}



\begin{theorem}[Theorem~1 (b) in \cite{e-s},
see also {\cite[Remark~3]{e-s}}]\label{thm:e-s}
The Smale invariant of $f$ is given by the formula:
$$\Omega(f) = \frac{1}{2}(3\sigma(V^4)+3t(\hat{f})
-3l(\hat{f})+L(f)).\quad\qed$$
\end{theorem}
\smallskip

As is stated in \cite[Remark~5]{e-s}, the above 
definition of 
$l$ makes sense for immersions of arbitrary closed 
oriented $3$-manifolds. 
Furthermore, the definition of $L$
can be extended to the case of framed immersions ---
immersions with trivial normal bundles and 
endowed with normal framings ---
of closed oriented $3$-manifolds, which are not 
necessarily $S^3$.
We review this extension here, since it will be necessary 
later.

\begin{definition}[{see \cite[Remark~5]{e-s}}]\label{dfn:Lnu}
Let $F$ be a framed immersion with a normal framing $\nu$, 
of a closed oriented $3$-manifold $M^3$ into $\R^5$.
Let $\Delta=\Delta(F)$ be the set of double points of $F$
and $F(\Delta)$ its image.
Let $n_1$ and $n_2$ be the two linearly independent normal 
vectors of $F$ determined by $\nu$.
Define the vector field $w$ along $F(\Delta)$ by 
$w(q) = n_1 (p_1) + n_1 (p_2)$, where 
$F(p_1)=F(p_2)=q$.
Then, define $L_\nu(F)$ to be the linking number in 
$\R^5$ of 
$F(M^3)$ and $F(\Delta)$ pushed off along $w$ out of 
$F(M^3)$.
\end{definition}

\begin{remark}
We will show later in the proof of Theorem~\ref{thm:B}
that $L_\nu(F)$ actually does not depend on 
the choice of the normal framing $\nu$
for an immersion  $F:M^3 \looparrowright \R^5$ 
with trivial normal bundle of 
any closed oriented $3$-manifold $M^3$ (see 
Remark~\ref{rmk:Lindep}).
\end{remark}

\vspace{5mm}  
\section{Immersions of $3$-manifolds into $\R^5$ --- 
Wu's result \\
and the geometry behind it}\label{sect:C}

In this section, we first recall the result of Wu 
\cite{wu} on 
the enumeration of regular homotopy classes of 
immersions of an oriented $3$-manifold 
$M^3$ into $\R^5$. 
Then we study the structure of $\Imm[M^3,\R^5]_0$ --- 
the set of 
regular homotopy classes of immersions with trivial normal 
bundles. We look especially into a certain second 
cohomology class which 
characterises (the regular homotopy class of) an immersion 
on (the regular neighbourhood of) the $2$-skeleton of 
$M^3$.


{\it We fix a trivialisation $\tau$ of $TM^3$ once and for 
all.}
Then the Smale-Hirsch theory provides a bijection
$\Imm[M^3, \R^5]\stackrel{\approx}{\to}[M^3, V_{5,3}]$.
Based on this bijection, Wu \cite{wu} has shown the 
following.

\begin{theorem}
[{Theorem~2 in \cite{wu}, see also \cite[Theorem~3]{li}}]
\label{thm:wu}
The normal Euler class $\chi_F$ for an immersion
$F:M^3 \looparrowright \R^5$ is of the form 
$2C$ for some $C\in H^2(M^3;\Z)$, and
for any $C\in H^2(M^3;\Z)$, there is an $F$ such that 
$\chi_F = 2C$.
Furthermore, 
$$\Imm[M^3, \R^5]_{\chi} \approx 
\coprod_{C\in H^2(M^3;\Z)\text{ with }2C=\chi}
H^3(M^3;\Z)/(4C\smile H^1(M^3;\Z)),$$
where $\Imm[M^3, \R^5]_{\chi}$ is the set of 
regular homotopy classes of immersions with 
normal Euler class $\chi\in H^2(M^3;\Z)$, 
{\rm ``}$\smile${\rm ''} represents the cup
product and {\rm ``}$\approx${\rm ''} denotes a 
bijection.\qed
\end{theorem}

\begin{remark}\label{rmk:C}
Let us consider the $S^1$-bundle $SO(5) \to V_{5,3}$ 
and choose the generator
$\Sigma^2 \in H^2(V_{5,3};\Z)$ for which the Euler class of 
this bundle 
coincides with $2 \Sigma^2$.
Let $\xi_F:M^3\to V_{5,3}$ be the map associated with
an immersion $F:M^3 \looparrowright \R^5$ (and the fixed 
trivialisation
$TM^3 = M^3 \times \R^3$).
The first part of the above theorem is based on the fact 
that the normal Euler class $\chi_F$ is the pullback 
$\xi_F^*(2\Sigma^2)$ for the generator $\Sigma^2$ for 
$H^2(V_{5,3};\Z) \approx \Z$. 
Hence, the class $C\in H^2(M^3;\Z)$ which appears in
the second part of the theorem is nothing but the pullback 
$\xi_F^*(\Sigma^2)$ of the generator $\Sigma^2$ of 
$H^2(V_{5,3};\Z) \approx \Z$. 
\end{remark}

Recall that $\Gamma_2(M^3)$ is 
the finite set $\{C\in H^2(M^3;\Z)\,|\,2C=0\}$.
Hereafter we shall concentrate on immersions with trivial 
normal bundles, i.e.\ on $\Imm[M^3, \R^5]_0$.
By Theorem~\ref{thm:wu} this set can be identified with 
$\Gamma_2(M^3) \times \Z$ (this identification depends 
on the trivialisation of $TM^3$).
Note that Theorem~\ref{thm:wu} can also be applied to 
the nonclosed $3$-manifold $\punc{M}^3$
and gives the following bijection:
$$\Imm[\punc{M}^3, \R^5]_0 \approx \Gamma_2(M^3).$$
Note that elements of order two in $H^2(M^3;\Z)$ 
can be naturally identified with those in 
$H^2(\punc{M}^3;\Z)$.

\begin{definition}
\label{dfn:C}
The projection $c: \Imm[M^3, \R^5]_0 \to \Gamma_2(M^3)$
is called {\it the Wu invariant\/} of the immersion of 
the {\it parallelised\/} $3$-manifold with trivial normal 
bundle. 
It can be described in two equivalent ways: 
\begin{itemize}
\item[(1)] $c(F) = \xi^*_F(\Sigma^2)$  or
\item[(2)] $c(F) = F|_{\punc{M}^3}$ if we use the 
identification 
$\Imm[\punc{M}^3, \R^5]_0 \approx \Gamma_2(M^3).$
\end{itemize}
\end{definition}

Next we shall give a more geometric description of the Wu 
invariant
$c(F) \in \Gamma_2(M^3)$.

A normal trivialisation $\nu$ of $F$ (together with the 
tangent trivialisation and the differential of $F$) 
defines a map 
$\pi_1(M^3) \to \pi_1(SO(5))$, i.e.\ an element  
$\tilde c_F$ in $H^1(M^3;\Z_2)$. If we change
$\nu$ by an element $z \in [M^3, SO(2)] = H^1(M^3;\Z)$, 
then the class $\tilde c_F$ changes by $\rho(z)$,
where $\rho$ is the mod $2$ reduction map $H^1(M^3;\Z) 
\to H^1(M^3;\Z_2)$.
Hence the coset of $\tilde c_F$ in 
$H^1(M^3; \Z_2)/\rho(H^1(M^3;\Z))$
does not depend on $\nu$.
The cokernel of $\rho$ can be identified with 
$\Gamma_2(M^3)$ 
by the canonical map induced by the Bockstein 
homomorphism (see below). 
Under this identification, the coset of $\tilde c_F$ 
corresponds
to the Wu invariant $c(F) \in \Gamma_2(M^3)=
\{C\in H^2(M;\Z)\,|\,2C = 0\}$.

Thus the Wu invariant describes the immersion on the 
$1$-skeleton
(and then actually also on the $2$-skeleton, since 
$\pi_2(SO(5)) = 0$),
if a trivialisation of the tangent bundle is fixed.
On the other hand for any immersion $M^3 \looparrowright 
\R^5$ there
is a trivialisation of the tangent bundle such that the 
Wu invariant of the
immersion under the given tangent bundle trivialisation 
is zero.
 
This geometric description will follow from the lemmas 
below.

Let $F:M^3\looparrowright \R^5$ be an immersion 
with trivial normal bundle. 
Choosing a normal framing $\nu$ for $F$, we 
naturally obtain 
a map $\varphi_{\nu,F}:M^3\to SO(5)$. 

\begin{lemma}\label{lem:Cx}
Let $x$ be the generator of $H^2(SO(5);\Z) \approx \Z_2$
and ${\varphi_{\nu,F}}^*$ the induced homomorphism 
${\varphi_{\nu,F}}^*:H^2(SO(5);\Z)\to H^2(M^3;\Z)$.
Then ${\varphi_{\nu,F}}^*(x)$ does not depend on the 
choice of 
the normal framing $\nu$, and is equal to $c(F)$.
\end{lemma}

\begin{proof} From the Gysin exact sequence of the 
$SO(2)$-bundle 
$p:SO(5) \to V_{5,3}$, it is easy to see that 
$p^*(\Sigma^2) = x$.
Hence ${\varphi_{\nu,F}}^*(x)$ 
is equal to $c(F)$ and does not depend on 
the choice of the normal framing $\nu$.
\end{proof}

Consider the following commutative diagram:
\[\begin{CD}
\! H^1(SO(5); \Z) \! @> \! \rho \! >> 
\! H^1(SO(5); \Z_2) \!
@> \! \beta \! >> \! H^2(SO(5); \Z) \!
@> \! {\times 2}\! >> \! H^2(SO(5); \Z) \! \\
@VV{\varphi_{\nu,F}}^*V @VV{\varphi_{\nu,F}}^*V @VV
{\varphi_{\nu,F}}^*V @VV{\varphi_{\nu,F}}^*V \\
\! H^1(M^3; \Z) \! @> \! \rho \! >> \! H^1(M^3; \Z_2) \!
@> \! \beta \! >> \! H^2(M^3; \Z) \! @>
\! {\times 2} \! >> \! H^2(M^3; \Z),\! \\
\end{CD}\] 
where each horizontal line is a part of 
the cohomology exact sequence associated with
the coefficient exact sequence
$0 \longrightarrow \Z \stackrel{\times 2}{\longrightarrow}
\Z \longrightarrow \Z_2 \longrightarrow 0$,
each vertical line is the homomorphism induced 
by $\varphi_{\nu,F}$, and each $\beta$ is 
the Bockstein homomorphism.

\begin{lemma}\label{lem:Cy}
Let $y$ be the generator of $H^1(SO(5);\Z_2) \approx \Z_2$.
Then $\beta({\varphi_{\nu,F}}^*(y))$ does not depend on 
the choice of 
the normal framing $\nu$, and is equal to $c(F)$.
\end{lemma}

\begin{proof}
Since $H^1(SO(5);\Z)=0$ and $H^2(SO(5);\Z) \approx \Z_2$
in the diagram above, we see that $\beta(y)=x$.
Hence ${\varphi_{\nu,F}}^*(x)=
\beta({\varphi_{\nu,F}}^*(y))$. Then
the result follows directly from Lemma~\ref{lem:Cx}.
\end{proof}

\begin{remark}
Lemma~\ref{lem:Cy} implies that the coset of the class
${\varphi_{\nu,F}}^*(y)\in H^1(M^3; \Z_2)$ 
modulo $\rho(H^1(M^3; \Z))$
does not depend on 
the choice of the normal framing $\nu.$ 

Note also that $$\Gamma_2(M^3) = \Ker 
(H^2(M^3; \Z)\stackrel{\times 2}{\longrightarrow}H^2(M^3; 
\Z))$$
is isomorphic to $H^1(M^3;\Z_2)/\rho(H^1(M^3;\Z))$ 
by the exact sequence above.
Thus, under this isomorphism, we have  
$c(F)={\varphi_{\nu,F}}^*(x)=[{\varphi_{\nu,F}}^*(y)]$, 
where the bracket means the coset
modulo the image of $\rho.$

Since the class $\tilde c_F \in H^1(M^3;\Z_2)$ mentioned 
above is nothing 
but the homomorphism 
$(\varphi_{\nu,F})_*:H_1(M^3;\Z)\to H_1(SO(5);\Z)=\Z_2$,
it is clear that 
\[{\varphi_{\nu,F}}^*(y)=y\circ (\varphi_{\nu,F})_*
=\tilde c_F:H_1(M^3;\Z)\to\Z_2,\]
where we consider the second $y$ as the identity 
homomorphism 
$H_1(SO(5);\Z)\to\Z_2$.
Thus, we have also the description $c(F)=[\tilde c_F]$.

Henceforth we often identify $\Gamma_2(M^3)$ with 
$H^1(M^3;\Z_2)/\rho(H^1(M^3;\Z))$ and 
consider $c(F)$ as an element of 
$H^1(M^3;\Z_2)/\rho(H^1(M^3;\Z))$.
\end{remark}

\begin{remark}\label{rmk:spin}
Since we have fixed a trivialisation of $TM^3$, the set
$H^1(M^3; \Z_2)$ can be identified with the set of 
homotopy classes
of trivialisations of $TM^3|_{\sk{2}{M}^3}$, 
where $\sk{2}{M}^3$ denotes the $2$-skeleton of $M^3$, 
i.e.\ with the
set of {\it spin structures\/} of $M^3$.
The above observations show that an immersion of $M^3$
into $\R^5$ with trivial normal bundle 
determines a spin structure
of $M^3$ up to elements of $\rho(H^1(M^3; \Z))$.
\end{remark}

Now we can state the main result of this section,
which claims that each $\Z$-component in 
$\Imm[M^3,\R^5]_0 \approx
\Z\amalg\Z\amalg\cdots\amalg\Z = \Gamma_2(M^3) \times \Z$ 
contains an embedding.

\begin{theorem}\label{thm:C}
For every element 
$C \in H^1(M^3; \Z_2)/\rho(H^1(M^3; \Z))$,
there exists an embedding $F:M^3\hookrightarrow\R^5$  
with the Wu invariant $c(F)=C$.
\end{theorem}

\begin{proof}
(A) First assume that $C=[0]$, where 
$0\in H^1(M^3; \Z_2)$.
Consider the spin structure of $M^3$ corresponding to 
$0\in H^1(M^3; \Z_2)$ 
under the identification in Remark~\ref{rmk:spin}.
Then by a result of Kaplan \cite{kaplan}, there
exists a framed 4-manifold $W$ such that
\begin{itemize}
\item[(1)] $\partial W = M^3$,
\item[(2)] the framing of $W$ restricted to 
the $2$-skeleton $\sk{2}{M}^3$ of $M^3$ coincides with
the given spin structure of $M^3$,
\item[(3)] $W$ has a special handlebody decomposition, 
consisting
of one 0-handle and some 2-handles attached
to the 0-handle simultaneously.
\end{itemize}
Since $W$ is a nonclosed spin 4-manifold, it is 
parallelisable 
and so it can be immersed into $\R^5$.
Furthermore, since $W$ has a $2$-complex as its spine, 
its immersion into $\R^5$ can be deformed into an 
embedding.
Take an embedding $\tilde{F}:W\hookrightarrow \R^5$  
and set $F=\tilde{F}|_{M^3}:M^3\hookrightarrow \R^5$.   

Now it suffices to show that $c(F)=[0]$.
Taking the normal framing $\nu$ for $F$ 
corresponding to the normal vector field of 
$M^3\subset W$, we see that the map 
$\varphi_{\nu,F}:M^3\to SO(5)$
restricted to the $2$-skeleton can be written as the 
following composition:
$$\varphi_{\nu,F}|_{\sk{2}{M}^3}:\,\sk{2}{M}^3\to W
\to SO(5),$$
where 
$W\to SO(5)$ is the map determined by 
the differential of $\tilde{F}$ using the given 
framing on $W$.
Thus  
${\varphi_{\nu,F}}^*(y)=0\in H^1(M^3;\Z_2)\approx 
H^1(\sk{2}{M}^3;\Z_2)$, 
since $H^1(W;\Z_2)=0$.

(B) In the case when $C$ cannot be written as 
$[0]$ with $0\in H^1(M^3;\Z_2)$, put $C=[\gamma]$ with 
$\gamma\neq 0\in H^1(M^3;\Z_2)$.

As in (A), consider the spin structure of $M^3$
corresponding to $\gamma\in H^1(M^3; \Z_2)$, and 
take an embedding $\tilde{G}:W\hookrightarrow\R^5$ 
of a $4$-manifold $W$ with $\partial W=M^3$
such that $W$ has the same properties as in (A).
Set $G=\tilde{G}|_{M^3}$.


In order to show that $C(G)=[\gamma]$, 
recall that we have fixed a trivialisation 
$\tau$ of $TM^3$ to obtain the map 
$c:\Imm[M^3,\R^5]_0\to\Gamma_2(M^3).$

If we take another trivialisation $\tau'$ of $TM^3$, then
we obtain another map 
$$c':\Imm[M^3,\R^5]_0\to\Gamma_2(M^3)$$
and with respect to this new trivialisation $\tau'$, 
we have a new map 
$$\varphi'_{\nu',G}:M^3\longrightarrow SO(5),$$
and a new induced homomorphism 
$${{\varphi'}_{\nu',G}}^*:H^1(SO(5);\Z_2)\longrightarrow 
H^1(M^3;\Z_2)$$
which satisfies  Lemma~\ref{lem:Cy}.
Furthermore, if we choose $\tau'$ so that 
$\tau'|_{\sk{2}{M}^3}$ corresponds to $-\gamma = \gamma 
\in H^1(M^3;\Z_2)$
as a spin structure, then we have 
\[\begin{array}{lcr}
{{\varphi'}_{\nu',G}}^*(y)={\varphi_{\nu',G}}^*(y) 
- \gamma 
&\text{ in }& H^1(M^3;\Z_2), 
\end{array}\]
and hence
\[\begin{array}{lcr}
c'(G)=c(G)-C&\text{ in }& H^1(M^3;\Z_2)/\rho(H^1(M^3;\Z)).
\end{array}\]
Since the restriction of
${\varphi'}_{\nu',G}: M^3\to SO(5)$
to the 
$2$-skeleton $\sk{2}{M}^3$ can be decomposed into the 
same composition   
$${\varphi'}_{\nu',G}|_{\sk{2}{M}^3}:\sk{2}{M}^3\to W
\to SO(5),$$
as in (A), by the same reason as in (A) we have 
$c'(G)=c(G)-C=0$ and hence $c(G)=C$.
Thus $G$ is a required embedding.
\end{proof}

\vspace{5mm}  
\section{Modifying immersions on a disk}\label{sect:sharp}

In this section we define an effective action of the group 
$\Imm[S^3,\R^5] \approx \Z$ on the set $\Imm[M^3, \R^5]_0$ 
and 
show that its orbits coincide with the sets of regular 
homotopy classes having the same Wu invariant $c$.

Let $M^3$ be a closed oriented $3$-manifold.
Let $D^3$ be a $3$-disk,
which from now on we will often
identify with the northern hemisphere
of the $3$-sphere $S^3$.
Fix an inclusion $D^3 \subset M^3$, 
and put $\punc{M}^3 = M^3\setminus\Int D^3$.
Suppose $F_0 : M^3 \looparrowright \R^5$ is an
immersion with trivial normal bundle
such that $F_0|_{D^3}$ coincides 
with the standard embedding $S^3 \hookrightarrow \R^5$
restricted to the northern hemisphere.
For an immersion $f : S^3 \looparrowright \R^5$, we may 
assume 
that $f$ restricted to the southern hemisphere is standard. 
Then consider the map
$$\begin{array}{cccc}
\sharp_{F_0} : &\Imm[S^3, \R^5] &\longrightarrow &
\Imm[M^3, \R^5]\\ 
 & f &\longmapsto & F_0 \sharp f, \\
\end{array}$$ 
where $(F_0 \sharp f)|_{\punc{M}^3} = F_0|_{\punc{M}^3}$ 
and $(F_0 \sharp f)|_{D^3} = f|_{D^3}$. 
The normal bundle of $F_0$ is trivial and if $F_0$ is 
modified on $D^3$, 
then its normal bundle does not change. Furthermore, 
$c(F_0)$ also does not change under this operation, since 
the invariant $c$ is determined by the immersion 
restricted to 
a neighbourhood of the $2$-skeleton of $M^3$.
Therefore, if we define, for $C\in H^2(M^3;\Z)$ with 
$2C=0$, 
$$\Imm[M^3, \R^5]_0^{C}:=\{F\in\Imm[M^3, \R^5]_0\, |
\, c(F)=C\},$$
then we can in fact define the map 
$$\sharp_{F_0} : \Imm[S^3, \R^5] \longrightarrow 
\Imm[M^3, \R^5]_0^{c(F_0)}.$$ 

The following proposition is an analogue of 
Proposition~3.1 in \cite{takase}.

\begin{proposition}\label{prop:bij}
The map
$$\sharp_{F_0} : \Imm[S^3, \R^5] \longrightarrow 
\Imm[M^3, \R^5]_0^{c(F_0)}$$
is a bijection.
\end{proposition}

\begin{proof}
As it has been mentioned in the paragraph just 
after Remark~\ref{rmk:C}, we have
\[\Imm[\punc{M}^3, \R^5]_0^{c(F_0)}\approx 
H^3(\punc{M}^3;\Z)= {0}.\]
This means that $\sharp_{F_0}$ is surjective. 

Let us prove the injectivity.
For two immersions $f$ and $g : S^3 \looparrowright \R^5$, 
we have
$$F_0 \sharp f \sim _r F_0 \sharp g \ 
\Leftrightarrow\ \xi_{F_0 \sharp f} 
\simeq \xi_{F_0\sharp g}:M^3\to V_{5,3}.$$
Note that $\xi_{F_0 \sharp f}$ and $\xi_{F_0\sharp g}$ 
are homotopic 
on the $2$-skeleton of $M^3$. 

Let $\Delta ^3 _{\xi_{F_0 \sharp f},\xi_{F_0 \sharp g}}$ be 
the difference $3$-cochain between $\xi_{F_0 \sharp f}$ and 
$\xi_{F_0 \sharp g}$, that is, 
the 3-cochain which assigns $\Omega(f)-\Omega(g)\in 
\pi _3(V_{5,3})$ 
to $D^3$ considered as a $3$-cell, 
and $0 \in \pi _3 (V_{5,3})$ to the other 3-cells.
Then we have 
$$\xi_{F_0 \sharp f} \simeq \xi_{F_0\sharp g}:M^3\to 
V_{5,3} \Leftrightarrow\ \Delta ^3 
_{\xi_{F_0 \sharp f},\xi_{F_0 \sharp g}} 
\text{ is a coboundary},$$
using the fact that $2c(F_0)=0$ and the following lemma.

\begin{lemma}[{Theorem~8A in \cite{whitney}, 
see also \cite[Lemma~2.1]{takase}}]
Two maps $\xi$ and $\eta : M^3 \to V_{5,3}$ are homotopic 
if and only if
\begin{itemize}
\item[$(a)$] $\xi^*(\Sigma^2)=\eta^*(\Sigma^2) \in 
H^2(M^3;\Z)$, and 
\item[$(b)$] there exist a $1$-cocycle $X^1$ and a 
$2$-cochain 
$Y^2$ such that 
$$\Delta^3 _{\xi ,\eta} = 4X^1 \smile \xi^*(\Sigma^2) 
+ \delta Y^2.\quad\qed$$
\end{itemize}
\end{lemma}
\smallskip

Note that $\Sigma^2$ is the generator of $H^2(V_{5,3};\Z)
\approx \Z$ and that 
$$\xi_{F_0 \sharp f}^*(\Sigma^2)=
\xi_{F_0 \sharp g}^*(\Sigma^2)
=\xi_{F_0}^*(\Sigma^2)=c(F_0)$$ 
(see Remark~\ref{rmk:C} and Definition~\ref{dfn:C}).

Now if $\Omega(f) = \Omega(g)$, then 
$\Delta^3 _{\xi_{F_0 \sharp f},\xi_{F_0 \sharp g}}$ 
(which we abbreviate here to $\Delta^3$) is trivial, 
hence a coboundary.
Conversely, if $\Delta^3$ is a coboundary, then 
there exists 
a $2$-cochain $Y^2$ such that $\delta Y^2=\Delta^3$. 
Thus, we have
$$Y^2(\partial D^3)=\delta Y^2(D^3)=\Delta^3(D^3)=
\Omega(f)-\Omega(g).$$
Furthermore, 
$$Y^2(\partial D^3)=-Y^2(\partial 
(M^3\setminus\Int{D^3}))
=-\Delta^3(M^3\setminus\Int{D^3})=0,$$
since $\Delta^3$ takes $0$ on all $3$-cells of 
$M^3\setminus\Int{D^3}$.
Therefore, we have 
\[\begin{split}
&\Delta^3 _{\xi_{F_0 \sharp f},\xi_{F_0 \sharp g}}
\text{ is a coboundary}\\
\Leftrightarrow\ &\Omega(f) = \Omega(g) \in  \pi _3 
(V_{5,3})\\
\Leftrightarrow\ &f \sim _r g :S^3 \looparrowright \R^5.
\end{split}\]
This completes the proof of Proposition~\ref{prop:bij}. 
(Note that this last argument
also shows that the map $\sharp_{F_0}$ is indeed 
well-defined.)
\end{proof}

\begin{remark}\label{rmk:act}
The proposition above implies that for each 
$C\in\Gamma_2(M^3)$,
the action of $\Imm[S^3, \R^5]$ on $\Imm[M^3, \R^5]_0^C$
\[\begin{array}{ccc}
\Imm[M^3, \R^5]_0^C\times\Imm[S^3, \R^5]&\to&
\Imm[M^3, \R^5]_0^C\\
(F, f)&\mapsto&F\sharp f
\end{array}\]
is effective and transitive.
\end{remark}

\vspace{5mm}  
\section{Geometric formulae for $\Imm[M^3,\R^5]_0$}
\label{sect:main}

In this section, we shall define an integer invariant 
$$i: \Imm[M^3, \R^5]_0 \to \Z,$$
which together with the Wu invariant $c$ will give a 
bijection 
$$(c,i): \Imm[M^3, \R^5]_0 \to \Gamma_2(M^3)\times \Z.$$  

Actually we shall give two different geometric expressions
for $i$ --- denoted by $i_a$ and $i_b$ --- which 
are analogues of the formulae in \cite{e-s}.
Then we show that they coincide.


\begin{definition}\label{dfn:alpha}
Let $M^3$ be a closed oriented 3-manifold.
We denote by $\alpha = \alpha(M^3)$ the
dimension of the $\Z_2$ vector space
$\tau H_1(M^3; \Z) \otimes \Z_2$, where
$\tau H_1(M^3; \Z)$ is the torsion subgroup of
$H_1(M^3; \Z)$. Note that the number of
elements in the set $\Gamma_2(M^3)$ is
equal to $2^\alpha$.
\end{definition}

The following lemma can be found
in \cite[Theorem~2.6]{k-m}. For reader's convenience we 
give a short proof here.

\begin{lemma}\label{lem:alpha}
Let $M^3$ be a closed oriented $3$-manifold and
$W$ a compact oriented spin $4$-manifold
with boundary $M^3$. Then the signature $\sigma(W)$
of $W$ has the same parity as $\alpha(M^3)$.
\end{lemma}

\begin{proof}
By \cite{kaplan}, there exists a 1-connected compact 
oriented
spin 4-manifold $V$ with $\partial V = M$ such that
the spin structures on $M^3$ induced by
$W$ and $V$ coincide with each other. Then by
Rohlin's theorem, $\sigma(V) - \sigma(W) = 
\sigma(V \cup -W)
\equiv 0 \pmod{16}$. Thus we have only to show the 
assertion for $V$
instead of $W$.

Since $V$ is 1-connected, we have the exact sequence
$$0 \longrightarrow H_2(M^3; \Z) \longrightarrow 
H_2(V; \Z) \stackrel{Q}{\longrightarrow} H_2(V, 
\partial V; \Z)
\longrightarrow H_1(M^3; \Z) \longrightarrow 0,$$
where the map $Q$ is identified with the intersection 
form
of $V$ through the Poincar\'e-Lefschetz duality $H_2(V, 
\partial V; \Z)
\approx H^2(V; \Z)$. By taking an appropriate
basis for $H_2(V; \Z)$, we may assume that the matrix
representative of $Q$ is of the form $Q_0 \oplus Q_1$,
where $Q_0$ is the zero form and $Q_1$ is nonsingular.
Then $Q_1$ can be regarded as a presentation matrix
of $\tau H_1(M^3; \Z)$ and its size has the same
parity as $\sigma(V)$. Now, since $V$ is spin,
its intersection form is of even type. Thus all
the diagonal entries of $(Q_1)_2$ are zero, where
$(Q_1)_2$ denotes the reduction modulo two. Now,
it is an easy exercise to show that
the size of such a matrix has the same parity
as the dimension of the $\Z_2$ vector space $\tau H_1(M^3;
\Z) \otimes \Z_2$ which it
presents. Thus the result follows.
\end{proof}

\begin{definition}
\label{dfn:ia}
Let $F : M^3 \looparrowright \R^5$ be 
an immersion with trivial normal bundle.
Let $W^4$ be any compact oriented $4$-manifold 
with $\partial W^4 = M^3$
and $\tilde{F}: W^4 \to \R^5$ a generic map
nonsingular near the boundary such that 
$\tilde F|_{\partial W^4} = F$. (We can choose such 
a generic map $\tilde{F}$,
since $F$ is an immersion with trivial normal bundle.)  
Denote the algebraic number of cusps of $\tilde{F}$ by 
$\#\Sigma^{1,1}(\tilde{F})$. Then define
$$i_a(F) = \frac{3}{2}(\sigma(W^4) - \alpha(M^3)) +
\frac{1}{2}\#\Sigma^{1,1}(\tilde{F}).$$
We will see later that this is always an integer.
\end{definition}

\begin{definition}
\label{dfn:ib}
Let $F : M^3 \looparrowright \R^5$ be 
an immersion with trivial normal bundle.
Let $W^4$ be any compact oriented $4$-manifold with 
$\partial W^4 = M^3$
and $\hat{F}: W^4 \to \R^6_+$ a generic map nonsingular 
near the boundary such that $\hat{F}^{-1}(\R^5) = \partial 
W^4$ and
$\hat{F}|_{\partial W^4} = F$. Then define 
$$i_b(F) = \frac{3}{2}(\sigma(W^4) - \alpha(M^3)) +
\frac{1}{2}(3t(\hat{F})-3l(\hat{F})+L_{\nu}(F)),$$
where $t(\hat{F})$ is the algebraic number of triple 
points of $\hat{F}$, 
$l(\hat{F})$ measures the linking of the singularity 
set of $\hat{F}$ 
with the rest of the image $\hat{F}(W^4)$ (see 
Definition~\ref{dfn:l}), 
and $L_{\nu}(F)$ is the linking in $\R^5$ 
of the image $F(M^3)$ and the double point set of $F$ 
pushed out of $F(M^3)$ with respect to a normal framing 
$\nu$ for $F$ 
(see Definition~\ref{dfn:Lnu}). It will be shown that
$i_b$ also takes only integer values.
\end{definition}

\begin{lemma}\label{lem:smootha}
The map $i_a$ is well-defined. That is, for an immersion 
$F : M^3 \looparrowright \R^5$
with trivial normal bundle, $i_a(F)$ does not depend
on the $4$-manifold $W^4$ or 
on the choice of the generic map $\tilde{F}: W^4 \to \R^5$.
\end{lemma}

\begin{proof}
Let $\tilde{F} : W \to \R^5$ and 
$\tilde{F}' : W' \to \R^5$ be generic maps
of the compact oriented $4$-manifolds $W$ and $W'$ 
respectively 
as in Definition~\ref{dfn:ia}.

First, consider the case where the normal vector fields 
along $F(M^3)$
in $\tilde{F}(W)$ and in $\tilde{F}'(W')$ are homotopic.  
Then after a suitable deformation near $F(M^3)$, we may 
regard
$\tilde{F}(W)\cup\tilde{F}'(W')$ as the image of 
a generic map of the closed $4$-manifold 
$W\cup -W'$ which is an immersion on a neighbourhood 
of $M^3$ in $W\cup -W'$.
Now we obtain  
$$\#\Sigma^{1,1}(\tilde{F}\cup -\tilde{F}')
+ 3\sigma(W\cup -W')=0$$
by Lemma~3 in \cite{szucs2}, which claims that 
for any generic map $g:X^4\to\R^5$ of a closed oriented 
$4$-manifold $X^4$, the equality  $\#\Sigma^{1,1}(g) 
+ 3\sigma(X^4) =0$ holds. 
Thus the result follows.

Suppose that the normal vector fields along $F(M^3)$
in $\tilde{F}(W)$ and in $\tilde{F}'(W')$ are not 
homotopic. 
We may assume that both normal vector fields are of 
equal constant length.
Let $M_+$ and $M_+'$ be formed by 
the endpoints of vectors of these two normal 
vector fields.
Then we may assume that $M_+$ and $M_+'$ 
intersect each other along a surface $S_+$ which is
the image of a surface $S\subset M^3$ by the map
associating to $x\in M^3$ the endpoint of the normal 
vector at $F(x)$.
Thus, 
$\tilde{F}(W)\cup \tilde{F}'(W')$ can be deformed, 
by a small modification 
near $F(M^3)$, into the image of a smooth generic 
map which has   
only Whitney umbrella singular points on $S$ and 
no other singularities nearby.
Now we can apply Lemma~3 in \cite{szucs2} again.
This completes the proof. 
\end{proof}

\begin{theorem}\label{thm:A}
Let $F$ and $G : M^3 \looparrowright \R^5$ be 
two immersions with trivial normal bundles such that
$c(F)=c(G)$. Then
$F$ and $G$ are regularly homotopic if and only if 
$i_a(F)=i_a(G)$. 
\end{theorem}

\begin{proof}
Let us first prove that $i_a(F)=i_a(G)$ if $F$ and $G$ 
are regularly homotopic.

Suppose that $F$ and $G$ are regularly homotopic.
Let $H:M^3\times[0,1]\looparrowright\R^5\times[0,1]$ be 
the track of 
a generic regular homotopy $H_t:M^3\looparrowright\R^5$ 
between $F$ and $G$.
Then we can take a normal vector field $\nu_t$ of $H_t$, 
which determines normal vector fields 
$\nu_F$ and $\nu_G$ of $F:M^3\looparrowright\R^5$ 
and $G: M^3\looparrowright\R^5$ for $t = 0$ and $t = 1$ 
respectively.

(The argument below follows that of Wells in 
\cite[p.~288]{wells}.) 
There is a small positive $\varepsilon$ such that 
on the collars 
$M^3 \times [0, \varepsilon)$ and $M^3 \times 
(1-\varepsilon, 1]$ 
the maps $(x,t) \mapsto F(x) + t\cdot\nu_F(x)$ and 
$(x,t) \mapsto G(x) + (t-1)\cdot\nu_G(x)$, 
respectively, are immersions.
Let us denote these immersions by 
$h_{[0,\varepsilon)}$ and 
$h_{(1-\varepsilon,1]}$ respectively.

The tangent bundle of the cylinder $M^3 \times I$ is 
$T(M^3 \times I) = 
\pi^*(TM^3) \oplus \varepsilon^1$, where $I = [0, 1]$,
$\pi: M^3 \times I \to M^3$ is the projection, 
and $\varepsilon^1$ is the trivial line bundle. 
We may assume that the regular
homotopy $H_t$ is such that $H_t = F$ for 
$t < 2\varepsilon$ 
and $H_t = G$ for $t > 1 - 2\varepsilon$.
Now we define the bundle monomorphism $\Phi: 
T(M^3\times I) \to \R^5$
as $dH_t \oplus \nu_t$, i.e.\ 
for a tangent vector $(\pi^*(v),s) \in T_{(x,t)}
(M^3 \times I)$, where
$v \in TM^3$ and $s \in \varepsilon^1$, put 
$\Phi(\pi^*(v),s) = dH_t(v) + s \cdot \nu_t$. 

This map is monomorphic on each fibre indeed, since 
$\nu_t$ is normal to
$dH_t(TM^3)$.
Furthermore, it coincides with the differentials of 
the immersions 
$h_{[0,\varepsilon)}$ and $h_{(1-\varepsilon, 1]}$ 
over the collars.
Hence --- by the relative version of Hirsch's theorem 
--- there is an
immersion $H': M^3 \times I \looparrowright \R^5$ 
which coincides with the given 
immersions on the collars.

Let $\tilde{F}:W_F\to\R^5$ and $\tilde{G}:W_G\to\R^5$ be
generic maps as in Definition~\ref{dfn:ia}.
Now via the smoothing process as in the proof of 
Lemma~\ref{lem:smootha},  
we obtain a generic map 
$\tilde{F}\cup H'\cup \tilde{G}:W_F\cup (M^3\times[0,1])
\cup -W_G\to\R^5$ 
of the closed oriented $4$-manifold $W_F\cup 
(M^3\times[0,1])\cup -W_G$.
By Lemma~3 in \cite{szucs2} we obtain that $i_a(F) 
= i_a(G).$ 

Now let us show that $i_a$ is injective on the set 
$\Imm[M^3,\R^5]_0^C$
of regular homotopy classes of immersions having the 
same Wu invariant $C$.
Notice that the analogous map $i_a^{S^3}:
\Imm[S^3, \R^5] \to \Z$
has been shown to be an isomorphism in \cite{e-s} 
(see Theorem~\ref{thm2.1} in the present paper).

Let us consider the effective and transitive group 
action described in
the previous section restricted to just one orbit:
$$\Imm[M^3,\R^5]_0^C \times \Imm[S^3,\R^5] \to 
\Imm[M^3,\R^5]_0^C.$$
Choose an embedding $F_0$ in 
$\Imm[M^3,\R^5]_0^C$ (by Theorem~\ref{thm:C},
such an embedding exists).
Let us consider the map
\[\begin{array}{ccc}
\Imm[S^3, \R^5] &\to& \Imm[M^3, \R^5]_0^C\\ 
f &\mapsto& F_0 \sharp f.
\end{array}\]
This map is bijective 
by Proposition~\ref{prop:bij}.
Compose it with the map 
$$i_a - i_a(F_0):\Imm[M^3, \R^5]_0^C \to \Q.$$
The resulting composition 
$\Imm[S^3, \R^5] \to \Imm[M^3, \R^5]_0^C \to \Q$
coincides with the map $i_a^{S^3}$. 

Indeed, for any embedded Seifert surface $W_{F_0}$ 
for $F_0$ we have
$$i_a(F_0) = \frac{3}{2}(\sigma(W_{F_0}) - 
\alpha(M^3)),$$
and 
for a singular Seifert surface of $F_0 \sharp f$
one can take 
the boundary connected sum of the Seifert surface 
$W_{F_0}$ for $F_0$ with 
any singular Seifert surface for $f$.

Hence the map $i_a$ is injective.
\end{proof}

\begin{lemma}\label{lem:smoothb}
For an immersion $F : M^3 \looparrowright \R^5$
with trivial normal bundle and a normal framing 
$\nu$, 
$i_b(F)$ does not depend on the $4$-manifold $W^4$ or 
on the choice of the generic map $\hat{F}: W^4 \to 
\R^6_+$.
\end{lemma}

\begin{proof}
Let $\hat{F} : W \to \R^6_+$ and 
$\hat{F}' : W' \to \R^6_+$ be generic maps
of the compact oriented $4$-manifolds $W$ and $W'$ 
respectively 
as in Definition~\ref{dfn:ib}.
Then, we obtain a generic map 
$\hat{F}\cup \hat{F}':W\cup -W'\to\R^6$. From the fact 
that for any generic map
$g:X^4\to\R^6$ of a closed oriented $4$-manifold $X^4$,
the equality $\sigma(X^4)-l(g)+t(g)=0$ holds 
(see \cite[Lemma~4]{e-s}),
the result follows.
\end{proof}

The following theorem states that the invariant $i_b$ 
is well-defined, and 
it is a regular homotopy 
invariant. In particular  it implies that $i_b(F)$ does 
not depend 
on the choice of the normal framing $\nu$ for $F$, either.

\begin{theorem}\label{thm:B}
Let $F$ and $G : M^3 \looparrowright \R^5$ be 
two immersions with trivial normal bundles such that
$c(F)=c(G)$. Then
$F$ and $G$ are regularly homotopic if and only if 
$i_b(F)=i_b(G)$. 
\end{theorem}

\begin{proof}
Let $F:M^3\looparrowright\R^5$ be an immersion with 
trivial normal bundle.
By Theorem~\ref{thm:C}, there is an embedding 
$F_0:M^3 \hookrightarrow \R^5$ with $c(F_0)=c(F)$.
Furthermore, by Proposition~\ref{prop:bij},
there exists an immersion $f:S^3 \looparrowright \R^5$ 
--- unique up to regular homotopy ---
such that 
$F \sim _r F_0 \sharp f$.
Let $H:M^3 \times I \looparrowright \R^5 \times I$
be the track of a generic regular homotopy between
$F$ and $F_0 \sharp f$.
Let $t(H)$ be the algebraic number of triple points of 
$H$.

Take a compact $4$-manifold $W_F$ and a generic map 
$\hat{F}:W_F\to\R^6_+$ 
as in Definition~\ref{dfn:ib} for $F$, and let
$\hat{F_0}:W_{F_0}\hookrightarrow \R^5$ be a
Seifert surface for $F_0$.
Then, we can deform
$\hat{F}(W_F)\cup H(M^3 \times I)\cup\hat{F_0}(W_{F_0})$ 
into the image of 
a generic map of the manifold $W_F \cup (M^3 \times I) 
\cup -W_{F_0}$ 
into $\R^5 \times I\subset \R^6_+$ bounded by 
$f:S^3 \looparrowright \R^5$ (see Figure~\ref{fig:1}).
Therefore, by Theorem~\ref{thm:e-s}, we have
\[\begin{split}
\Omega (f) &= \frac{1}{2}(3\sigma(W_F\cup (M^3\times I)
\cup -W_{F_0}) -3t(H) +3t(\hat{F}) -3l(\hat{F}) + L(f))\\
           &= {\rule[0cm]{0cm}{0.7cm}\frac{3}{2}}
(\sigma (W_F) -\sigma (W_{F_0}) +t(\hat{F}) -l(\hat{F})) 
+\frac{1}{2}(L(f) - 3t(H)).\\ 
\end{split}\]

\begin{figure}[htbp]
\begin{center}
 \includegraphics[width=\linewidth,
height=0.6\textheight,keepaspectratio]{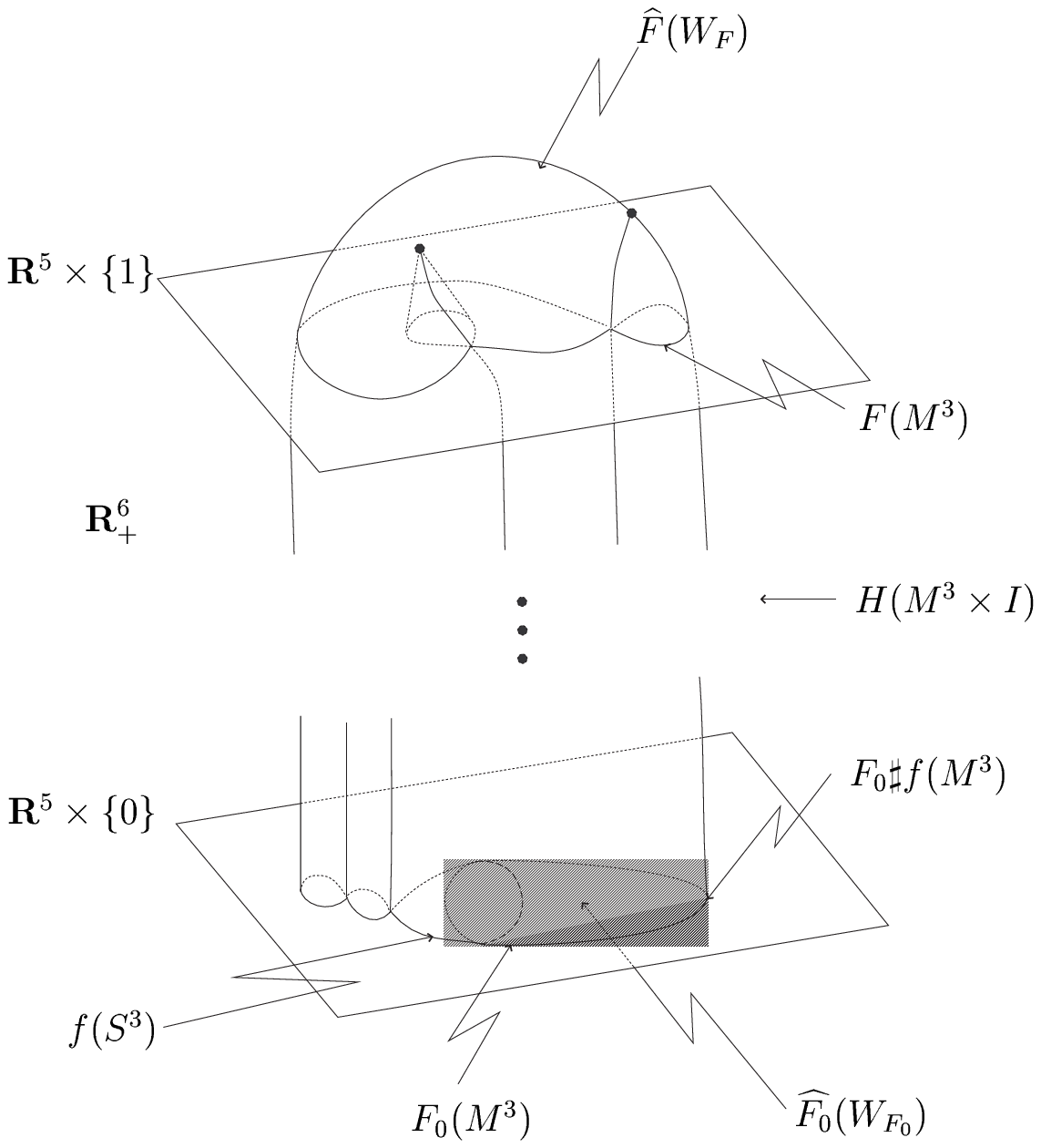}
\end{center}
\caption{}
\label{fig1}
\end{figure}

Here we need the following lemma.

\begin{lemma}\label{lem:L+3t}
We have
$$L(f)=L_\nu(F) +3t(H)$$
for any normal framing $\nu$ for $F$.
\end{lemma}

\begin{proof}
We can extend $\nu$ to the whole cylinder $M^3 \times I$ 
using the regular homotopy $H.$
We shall denote this extended normal field also by $\nu$.
Since on $D^3\times\{0\}\subset M^3\times\{0\}$,
the normal field $\nu$ coincides with the homotopically 
unique 
normal framing of $(F_0 \sharp f)|_{D^3}$, we have
$$L_{\nu}(F_0 \sharp f) = L(f).$$
(Indeed, recall that for an immersion the invariant 
$L_{\nu}$
was defined as the linking number of the image of the 
immersion with 
the curves formed by the double points of the immersion 
pushed off out of the image of the immersion into the 
direction defined by $\nu$.
Since $F_0$ is an embedding 
the equality follows.)

Furthermore, one can show that 
$$L_{\nu}(F_0 \sharp f) -L_{\nu}(F) = 3t(H)$$
in the same way as in \cite{ekholm1}.
This completes the proof of Lemma~\ref{lem:L+3t}.
\end{proof}

Thus, $L_\nu(F)$ does not depend on the choice of 
$\nu$ and we obtain 
\[\begin{split}
\Omega (f) &= \frac{3}{2}(\sigma (W_F) -\sigma (W_{F_0}) 
+t(\hat{F}) 
               -l(\hat{F})) +\frac{1}{2}L_{\nu}(f)\\   
           &= {\rule[0cm]{0cm}{0.7cm}\frac{1}{2}}(3\sigma 
(W_F) +3t(\hat{F}) 
               -3l(\hat{F}) +L_\nu(F)) - \frac{3}{2}\sigma 
(W_{F_0}).       
\end{split}\]
Since the regular homotopy class of $F|_{\punc{M}^3}$
is determined by $c(F)$, 
again by Proposition~\ref{prop:bij},
the regular homotopy class of $F$ is completely determined 
by $\Omega (f)$, hence by $i_b(F)=\Omega (f)+ 
3(\sigma (W_{F_0}) - \alpha(M^3))/2$. 
This completes the proof of Theorem~\ref{thm:B}.
\end{proof}

\begin{remark}\label{rmk:Lindep}
Theorem~\ref{thm:B} implies that $i_b(F)$ does not 
depend on the choice of 
$\hat{F}$ or on the normal framing $\nu$ for $F$. Hence 
$L_{\nu}(F)$ 
does not depend on the choice of the normal framing $\nu$ 
for $F$, either.
Thus, we can use the notation $L(F)$ for an immersion $F$
with trivial normal bundle of 
a closed oriented $3$-manifold from now on; then,
the formula displayed in 
Definition~\ref{dfn:ib} can be written as 
$$i_b(F) = 
\frac{3}{2}(\sigma(W^4) - \alpha(M^3)) +
\frac{1}{2}(3t(\hat{F})-3l(\hat{F})+L(F)).$$
\end{remark} 

\begin{remark}\label{rmk:coin}
We see easily that for an embedding $F_0:M^3 
\hookrightarrow\R^5$
we have $i_a(F_0)=i_b(F_0)$; furthermore, this is always
an integer by Lemma~\ref{lem:alpha}.
Therefore, by Theorems~\ref{thm:A}, \ref{thm:B}, 
\ref{thm2.1} and \ref{thm:e-s},
$i_a$ and $i_b$ take integer values for general immersions 
(with trivial normal bundles).  
Furthermore, considering the action of 
$\Imm[S^3,\R^5]$ on $\Imm[M^3,\R^5]_0$,
we see that in fact $i_a=i_b$ (using Theorems~\ref{thm2.1} 
and \ref{thm:e-s}).
Thus, we obtain the surjective map
$$i:\Imm[M^3,\R^5]_0\longrightarrow\Z$$
defined in two ways as $i_a$ and as $i_b$,
so that the map
$$(c,i):\Imm[M^3,\R^5]_0\longrightarrow\Gamma_2(M^3)
\times\Z$$
gives a bijection.
\end{remark}

\begin{corollary}\label{cor:emb}
Let $M^3$ be a closed oriented $3$-manifold such that 
$H^2(M^3;\Z)$ has no $2$-torsion. Then 
two embeddings $M^3 \hookrightarrow \R^5$
are regularly homotopic if and only if 
they have Seifert surfaces with the same signature.
\end{corollary}

\begin{proof}
Since $H^2(M^3;\Z)$ has no $2$-torsion, 
$c(F)$ is equal to $0\in H^2(M^3;\Z)$
for every embedding $F : M^3 \hookrightarrow \R^5$.
If we take a Seifert surface 
$W^4\hookrightarrow \R^5$ for an embedding $F$,
then clearly $i_a(F)=i_b(F)=3(\sigma(W^4) - 
\alpha(M^3))/2$.
Thus the result follows directly from Theorems~\ref{thm:A} 
and \ref{thm:B}.
\end{proof}

\begin{remark}
As a consequence of Theorems~\ref{thm:A} and \ref{thm:B},
we have that {\it if 
two embeddings $F$ and $G : M^3 \hookrightarrow \R^5$ of 
an arbitrary closed oriented $3$-manifold $M^3$
into $\R^5$ are regularly homotopic, then their Seifert 
surfaces 
have the same signature}. In fact, we can also prove this 
fact
by using a simpler argument as follows.

Let $\tilde{F}:W_F\hookrightarrow\R^5$ and 
$\tilde{G}:W_G\hookrightarrow\R^5$
be Seifert surfaces for $F$ and $G$ respectively.
Let $H:M^3\times I\looparrowright\R^5\times I$ be 
the track of a generic regular homotopy between $F$ 
and $G$.
Then, using $\tilde{F}$, $\tilde{G}$, and $H$, we 
obtain a generic
immersion $h:W_F\cup(M^3\times I)\cup -W_G
\looparrowright\R^6$ 
after an appropriate smoothing process. 
Put $X^4:=W_F\cup(M^3\times I)\cup -W_G$. 
Let $\chi \in H^2(X^4;\Z)$ be the normal Euler class of 
the immersion $h$ and $e \in H_2(X^4;\Z)$ its 
Poincar\'e dual.
Then, $-p_1 (X^4) = \chi \smile \chi$ and $-e$ is 
represented by 
the set of double points  
$\Delta:= \{x\in X^4\, |\,h(x)=h(y)
\text{ for some }y \neq x \}$ of $h$ 
(see \cite[p.~44, Lemmas~1 and 4]{kirby}).
Since $\Delta \subset M^3 \times I \subset X^4$, we have
\begin{center}
\begin{tabular}{ccl}
$-p_1 (X^4)$ &=& $\chi \smile \chi \ \ \in H^4(X^4;\Z)
                                          \approx \Z$\\
             &=& $e \cdot e \ \ \ \ \in H_0(X^4; \Z)
                                          \approx \Z$\\
             &=& the self-intersection of $\Delta$ in 
                                               $X^4$\\
             &=& the self-intersection of $\Delta$ in 
                                       $M^3 \times I$\\
             &=& $0$.
\end{tabular}\\
\end{center}
Thus, we have $\sigma(X^4)=-p_1(X^4)/3=0$ and 
hence $\sigma (W_F) = \sigma (W_G)$.
\end{remark}

\vspace{5mm}  
\section{Embeddings of $T^3$ into $\R^5$}
\label{sect:t3}

In this section, we prove the 
two surprising corollaries formulated in the
introduction. Recall that these corollaries claim 
roughly that for the 
$3$-dimensional torus $T^3$
``more regular homotopy classes of immersions into 
$\R^5$ contain embeddings 
than for the $3$-sphere''. (The proofs will show 
that there 
are twice as many regular homotopy classes
containing embeddings for the torus as for the sphere.)

Since $\Gamma_2(T^3) = 0$, 
the set $\Imm[T^3, \R^5]_0$ can be identified with 
the set of integers.  
Therefore, by Proposition~\ref{prop:bij}, we see that 
$$\sharp_{F_0} :\Imm[S^3, \R^5]\longrightarrow
\Imm[T^3, \R^5]_0$$
is a bijection for any immersion $F_0: T^3 
\looparrowright \R^5$
with trivial normal bundle. Furthermore,
$i:\Imm[T^3, \R^5]_0\to\Z$ in Section~\ref{sect:main} 
gives 
a complete invariant of regular homotopy.

Obviously there exists an embedding  
$F_0: T^3 \hookrightarrow \R^5$ 
which bounds a $4$-manifold of signature $0$ 
(e.g.\ $D^2 \times S^1 \times S^1$) in $\R^5$.
Furthermore, we show the following.

\begin{proposition}\label{prop:t3}
\begin{itemize}
\item[{\rm (a)}] There exists an embedding $F_8 :T^3 
\hookrightarrow \R^5$ 
having a Seifert surface $W_{F_8} \hookrightarrow 
\R^5$ with signature $8$. 
\item[{\rm (b)}] The signature of any Seifert surface of 
any embedding
$T^3 \hookrightarrow \R^5$ is divisible by $8$.
\end{itemize}
\end{proposition}

\begin{proof}
(a) It is known that there exists 
a spin $4$-manifold $W = W_{F_8}$ of signature $8$ with 
$\partial W = T^3$, see
\cite{kirby} for example.
Furthermore, $W$ can be chosen to have 
a special handlebody decomposition with one 
$0$-handle and some $2$-handles with even framings 
(see \cite{kaplan}).
As we have seen in the proof of Theorem~\ref{thm:C},
such a $4$-manifold $W$ embeds into $\R^5$.
Such an embedding,
restricted to the boundary $\partial W$,
gives a required embedding  $F_8: T^3 \hookrightarrow \R^5$.

(b) It is shown also in \cite{kirby} that $T^3$ with any 
spin structure
spin-bounds either the solid torus $D^2 \times S^1 \times 
S^1$ or the 
above mentioned $W_{F_8}$.
Now let $V$ be any Seifert surface
of any embedding $T^3 \hookrightarrow \R^5$.
Then $V$ induces a spin structure on $T^3$, which is 
induced also from a spin structure either on the solid 
torus or on $W_{F_8}$.

Therefore, either $V \cup (-D^2 \times S^1 \times S^1)$ or 
$V \cup -W_{F_8}$ carries a spin structure, and so by 
Rohlin's theorem 
its signature is divisible by $16$.
By Novikov's additivity these signatures
are $\sigma(V) - \sigma(D^2 \times S^1 \times S^1) = 
\sigma (V)$
and $\sigma(V) - \sigma(W_{F_8}) = \sigma (V) - 8$
respectively.
In both cases $\sigma(V)$ is divisible by $8$.
\end{proof}

\begin{corollary}\label{cor:t3sz}
There exist two embeddings $F_0$ and 
$F_8:T^3\hookrightarrow\R^5$ such that 
there is an immersion $h:S^3\looparrowright \R^5$ with 
$F_0 \sharp h \sim_r F_8$, but 
$h$ cannot be chosen from a regular homotopy class 
containing an embedding.
\end{corollary}

\begin{proof}
Let $F_0$ and $F_8$ be the embeddings described above
in Proposition~\ref{prop:t3};
i.e.\ they have Seifert surfaces%
$W_{F_0}$ and $W_{F_8}$ of signatures $\sigma(W_{F_0})=0$
and $\sigma(W_{F_8})=8$ respectively.

From the fact that%
$$\sharp_{F_0} :\Imm[S^3, \R^5]\longrightarrow
\Imm[T^3, \R^5]_0$$
is a bijection, there is an immersion $h:S^3
\looparrowright \R^5$ such that  
$F_0 \sharp h \sim_r F_8$.
Let us assume that $h$ is regularly homotopic to 
an embedding.
Then the following sequence of equalities hold:
$$
12 = \frac{3}{2}\sigma(W_{F_8}) = i(F_8) = 
i(F_0 \sharp h) =
i(F_0) + \Omega(h) = 0 + 24k
$$
for an integer $k \in \Z$, which is a contradiction.
Here the fourth equality 
$i(F_0 \sharp h) =
i(F_0) + \Omega(h)$
follows from the obvious remark that
for a Seifert surface of the connected sum 
$F_0 \sharp h$ one can choose
the boundary connected sum
of the Seifert surfaces of the embeddings $F_0$ and $h$.
The Smale invariant $\Omega (h)$ is a multiple of $24$
by the result of Hughes and Melvin \cite{h-m} recalled in 
Section~\ref{sect:rev}.
%
\end{proof}

\begin{corollary}\label{cor:t3}
There exists an immersion $h : S^3 \looparrowright \R^5$ 
not regularly homotopic to any embedding $S^3 
\hookrightarrow \R^5$ such that
for any embedding $E:T^3 \hookrightarrow \R^5$,
the connected sum $E \sharp h$ of $F$ and $h$ is again 
regularly homotopic to 
an embedding $T^3 \hookrightarrow \R^5$. 
\end{corollary}

\begin{proof}
Let $h : S^3 \looparrowright \R^5$ be any immersion 
with Smale invariant $\Omega(h)=12$. 
Then $h$ is not regularly homotopic to an embedding
again by the result of Hughes-Melvin \cite{h-m}.
We show that for any embedding $E: T^3 \hookrightarrow 
\R^5$ the connected sum $E \sharp h$ is 
regularly homotopic to an embedding.
For this purpose we first compute the $i$-invariant of 
$E \sharp h$
and then produce an embedding with the same $i$-invariant.
Since the $i$-invariant determines the regular homotopy 
class completely,
we obtain that $E \sharp h$ is regularly homtopic to the 
embedding.

The following sequence of equalities hold:
$$i(E \sharp h) = i(E) + \Omega (h) = \frac{3}{2} 
\sigma(W_E) + 12 =
\frac{3}{2}\cdot 8k + 12 = 12 (k+1).$$
Here the first equlaity holds, since we can form 
the boundary connected sum of the Seifert surface 
$W_E$ of $E$ with that
of $h$ even if the latter is a ``singular'' one 
(i.e.\ it is a generic map
with boundary $h$, not necessarily an embedding). 
The number
$k = \sigma (W_E)/8$ is an integer by part (b) of 
Proposition~\ref{prop:t3}.
Let $e_n: S^3 \hookrightarrow \R^5$ be an embedding 
with Smale invariant 
$\Omega (e_n) = 24 n$, where $n = k/2$ if $k$ 
is even and 
$n = (k+1)/2$ if $k$ is odd.

Now we see that if $k$ is even, then $i(F_8 \sharp e_n) 
= i(E \sharp h)$.
If $k$ is odd, then $i(F_0 \sharp e_n) = i(E \sharp h)$.
\end{proof}

\vspace{5mm}  
\section{Number of cusps of singular Seifert surfaces}
\label{sect:cusp}

The present section is not along the main line of the 
paper.
Here as a byproduct of our formulae for $i_a$ and $i_b$, 
their correctness,
integrality and coincidence, we prove a few simple 
corollaries on the 
number of cusps a singular Seifert surface can have.
%

\begin{corollary}\label{cor:cusp1}
In the notation of Definition~{\rm \ref{dfn:ia}}, 
the residue class of the 
algebraic number $\#\Sigma^{1,1}$ of cusps of $\tilde F$ 
modulo $3$ does not depend on the $4$-manifold $W^4$
or on the choice of the map $\tilde F$, and it depends 
only on the immersion
$F$ of the $3$-manifold $M^3$ into $\R^5$.\qed
\end{corollary}

Naturally arises the problem: how to read off this modulo 
$3$ residue class from $F$?
It turns out that $\# \Sigma^{1,1}(\tilde F)\equiv 
L_{\nu}(F)\pmod 3$
(see Definition~\ref{dfn:Lnu}). 
This follows from the definitions 
of the invariants $i_a$ and $i_b$ 
(Definitions~\ref{dfn:ia} and \ref{dfn:ib}) 
and their coincidence (see Remark~\ref{rmk:coin}).

\begin{corollary}\label{cor:cusp2}
Let $W^4$ be a closed connected spin $4$-manifold, and 
$f: W^4 \to \R^5$ a generic map.
It is well-known that the singularity set is a surface 
$\Sigma(f)$ with isolated cusps on it.
Then the number of cusps on each component of $\Sigma(f)$ 
is even.
\end{corollary}

\begin{proof}
Otherwise, restricting $f$ to the closure of an appropriate 
neighbourhood of a component of 
$\Sigma(f)$, we would get a contradiction to the 
integrality of the invariant 
$i_a$ and Lemma~\ref{lem:alpha}.
\end{proof}

\begin{corollary}\label{cor:cusp3}
In the conditions of the previous corollary, let 
$M^3$ be a 
null-homologous oriented submanifold of $W^4$ disjoint 
from $($the closure of\/$)$ the double point set of $f$. 
Then in both parts of the set $W\setminus M$ 
the algebraic number of cusps is divisible by $6$.\qed
\end{corollary}


\end{document}